\documentclass[12pt]{amsart}

\usepackage{hyperref}
\usepackage{amsthm}
\usepackage{amsmath}
\usepackage{amssymb}
\usepackage{bbm}
\usepackage{physics}
\usepackage{mathtools}
\usepackage{enumitem}

\makeatletter
\@namedef{subjclassname@2020}{%
  \textup{2020} Mathematics Subject Classification}
\makeatother

\newtheorem{thm}{Theorem}[section]
\newtheorem{cor}[thm]{Corollary}

\newtheorem{lem}[thm]{Lemma}

\begin{document}

\baselineskip=17pt

\title{Minkowski bases, Korkin-Zolotarev bases and Successive Minima}
\author{Shvo Regavim}
\address{School of Mathematical Sciences, Tel Aviv University, Tel Aviv 69978, Israel}
\email{shvoregavim@mail.tau.ac.il}
\date{\today}



\maketitle

\section{Introduction and Statement of Results}

A \emph{lattice} is a discrete subgroup in $\mathbb{R}^m$, and its \emph{rank} is the dimension of the vector space it spans. Every lattice of rank $n$ has a \emph{basis}, that is a sequence of $n$ vectors $b_1, b_2, \dots, b_n \in L$ that generate $L$ as a free abelian group. Equivalently, $b_1, \dots, b_n$ are linearly independent and
$$\mathrm{span}_{\mathbb{Z}} \left( b_1, \dots, b_n \right) = L.$$
The \emph{i-th successive minimum} of $L$, denoted by $\lambda_i (L)$ or $\lambda_i$ is the smallest real number such that there are $i$ linearly independent vectors in $L$ of length at most $\lambda_i (L)$.

A sequence $v_1, v_2, \dots, v_k$ of elements of $L$ is called \emph{primitive} if there exist $v_{k+1}, \dots, v_n \in L$ such that $v_1, \dots, v_n$ is a basis of $L$. It is not difficult to show that an equivalent criterion is
$$\mathrm{span}_{\mathbb{Z}} (v_1, \dots, v_k) = \mathrm{span}_{\mathbb{R}} (v_1, \dots, v_k) \cap L.$$
Let $e_1, \dots, e_m$ be the standard basis vectors in $\mathbb{R}^m$. Given a sub-vector space $V \subset \mathbb{R}^m$ its orthogonal complement in $\mathbb{R}^m$ will be denoted by $V^{\perp}$.

The problem of selecting for a given lattice a basis which satisfies "good" properties is called \emph{reduction theory}. Classically, the theory was motivated by the problem of the finding the minimum of a positive definite integral quadratic form. The theory was started by Lagrange \cite{Lagrange} for binary forms, and Hermite \cite{Hermite} started the study of reduction theory in arbitrary dimensions.

Many different notions of reduction have been studied, but we will describe two of these notions which will be relevant to us.

A basis $b_1, \dots, b_n$ is \emph{reduced in the sense of Minkowski} or a \emph{Minkowski basis} if $b_i$ is the shortest vector such that $b_1, \dots, b_i$ is a primitive tuple. This is essentially a greedy algorithm which tries to choose the shortest basis possible. The even greedier algorithm where at each step we choose the shortest vector $b_i$ such that $b_1, \dots, b_i$ are linearly independent unfortunately does not always work. For example, take the lattice
$$D_{5}^{*} = \mathrm{span}_{\mathbb{Z}} \left( e_1, e_2, e_3, e_4, \frac{e_1 + e_2 + e_3 + e_4 + e_5}{2} \right)$$
which is the dual of the root lattice $D_5$. All vectors in $D_{5}^{*}$ are of length at least 1, and $e_1, e_2, e_3, e_4, e_5$ is a sequence of linearly independent vectors of minimal length, however they do not form a basis as $\frac{e_1 + e_2 + e_3 + e_4 + e_5}{2}$ is not an integral linear combination of them.

A basis $b_1, \dots, b_n$ is \emph{reduced in the sense of Korkin-Zolotarev} or a \emph{Korkin-Zolotarev basis} if for all $i, \ b_i$ is the shortest vector in $L$ among those that minimize $\norm{\pi (b_i)}$ among vectors $b_i \in L, \ b_i \not\in \mathrm{span}_{\mathbb{R}} \left( b_1, \dots, b_{i - 1} \right)$, where $\pi$ is the orthogonal projection on $\mathrm{span}_{\mathbb{R}} \left( b_1, \dots, b_{i - 1} \right)^{\perp}$.

In both the Minkowski reduction and the Korkin-Zolotarev reduction there is the question of what to do in the case of ties, and indeed this will occur in the examples we will give later. The answer is just to choose one of the vectors at random. However, at least in our context we may just perturb the lattice slightly to ensure that the only tie possible is between $v$ and $- v$ for some vector $v$, and then the Minkowski and Korkin-Zolotarev bases are unique up to signs.

From now on we will always let $v_1, \dots, v_n$ and $u_1, \dots, u_n$ denote the Minkowski and Korkin-Zolotarev bases respectively. Van der Waerden \cite{Van Der Waerden} showed that
\begin{align*}
    \norm{v_i}^2 \leq \Delta_i \lambda_{i}^2 \ \mathrm{with} \ \Delta_i = \max \left\{ 1, \left( \frac{5}{4} \right)^{i - 4} \right\}.
\end{align*}
In \cite[~p. 22-23]{Schurmann} Sch\"{u}rmann conjectures that for $i \geq 4$
$$\norm{v_i}^2 \leq \frac{i}{4} \lambda_{i}^2$$
which is a significant strengthening of van der Waerden's bound, replacing an exponential factor with a linear one. The lattice
\begin{align*}
    D_{n}^{*} = \mathrm{span}_{\mathbb{Z}} \left( e_1, \dots, e_{n - 1}, \frac{e_1 + \cdots + e_n}{2} \right)
\end{align*}
which is the dual of the root lattice $D_n$ shows that the above conjecture, if correct, is tight.

We will give some improvements on van der Waerden's bound for $v_k$, which were previously the best known bounds. In particular we show that Sch\"{u}rmann's conjecture mentioned above is true for $i = 6, 7$, that is

\begin{thm}\label{Bound}

For any lattice $L$, and for $k = 6, 7$ we have
$$\norm{v_k}^2 \leq \frac{k}{4} \lambda_{k}^2.$$
Furthermore, equality occurs if and only if
$$\mathrm{span}_{\mathbb{Z}} \left( v_1, \dots, v_k \right)$$
is similar to $D_{k}^{*}$ for $k = 6, 7$ respectively.

\end{thm}

Martinet \cite{Martinet} has shown that for $n \leq 8$, any basis $w_1, w_2, \dots, w_n$ of $L$ satisfies the inequality
$$\prod_{k = 1}^{n} \frac{\norm{w_k}^2}{\lambda_{k}^2} \leq \frac{n}{4}$$
with equality possible only when $L$ is similar to $D_{n}^{*}$. In particular, this shows that for $n \leq 8$
$$\norm{v_n}^2 \leq \frac{n}{4} \lambda_{n}^2.$$
Notice that this does not directly imply Theorem \ref{Bound}, as in our theorem there is no dependence on the dimension. The approach in \cite{Martinet} is very different from ours, and involves deformation arguments as well as Watson's index theory.

As we describe in the proof of Theorem \ref{Bound}, from this we can deduce a minor improvement to the bound in general:

\begin{cor}\label{Bound2}

For all $k \geq 8$ we have
$$\norm{v_k}^2 \leq \frac{608}{625} \left( \frac{5}{4} \right)^{k - 4} \lambda_{k}^2.$$
$\frac{608}{625} = 0.9728$.
\end{cor}

It seems that to get a significantly better bound a more refined technique is required.

In the other direction, since by definition $\norm{v_1} \leq \norm{v_2} \leq \cdots \leq \norm{v_n}$ and $v_1, \dots, v_i$ are linearly independent, we have
$$\norm{v_i}^2 \geq \lambda_{i}^2.$$
In fact, we have an even better lower bound. Let $\overline{\lambda_i} (L)$ or $\overline{\lambda_i}$ be the smallest real number such that there exist a primitive sequence of vectors $b_1, \dots, b_i$ all of which are of length at most $\overline{\lambda_i}$. As $v_1, \dots, v_i$ are primitive, then
$$\norm{v_i}^2 \geq \overline{\lambda_i}^2.$$
Clearly $\overline{\lambda_i} \geq \lambda_i$, so this is a better bound than before, and equality does not necessarily hold, as the example of $D_{5}^{*}$ shows where $\lambda_5 = 1, \norm{v_5} = \overline{\lambda_5} = \frac{5}{4}$. In a recent course on Geometry of Numbers \cite{Geometry} the question was raised to find a lattice $L$ such that $\norm{v_n} > \overline{\lambda_n}$. Given a basis $b_1, \dots, b_n$ of $L$, let us call it a \emph{shortest basis} of $L$ if
$$\max \left\{ \norm{b_1}, \dots, \norm{b_n} \right\} = \overline{\lambda_n}$$
that is the basis $b_1, \dots, b_n$ minimizes the maximum length of an element in the basis. In this terminology, the question is to find a lattice $L$ where the Minkowski basis is not a shortest basis of $L$. We will show an even stronger statement:

\begin{thm}\label{4b}

There exists a lattice $L$ such that
$$\norm{v_n} > \max \left\{ \norm{u_1}, \dots, \norm{u_n} \right\} = \overline{\lambda_n}.$$

\end{thm}

This answers a question of Sch\"{u}rmann (personal communication) of whether the longest vector in the Minkowski reduction can be longer than the longest vector in the Korkin-Zolotarev reduction. 

Sch\"{u}rmann \cite[~p. 23]{Schurmann} compares the conjecture $\norm{v_k}^2 \leq \frac{k}{4} \lambda_{k}^2$ to the corresponding known bounds on the Korkin-Zolotarev basis
$$\frac{4}{i + 3} \lambda_{i}^2 \leq \norm{u_i}^2 \leq \frac{i + 3}{4} \lambda_{i}^2$$
which were proven by Lagarias, Lenstra and Schnorr \cite{Bases}. Of course, if one could show an inequality of the form
$$\norm{v_i}^2 \leq C \max \left\{ \norm{u_1}^2, \dots, \norm{u_i}^2 \right\}$$
for some constant $C$ then this would show that $\norm{v_i}^2 \leq \frac{C (i + 3)}{4} \lambda_{i}^2$ which is up to a constant the conjectured asymptotic behaviour. However, this is not the case. In fact,

\begin{thm}\label{4b+}

There exist a sequence of lattices $\left( L_i \right)_{i = 1}^{\infty}$ of ranks $a_i \uparrow \infty$ such that
$$\norm{v_{a_i} (L_i)}^2 \geq \frac{c \cdot a_{i}^{\frac{1}{3}}}{\log^{\frac{2}{3}} a_i} \max \left\{ \norm{u_1 (L_i)}^2, \dots, \norm{u_{a_i} (L_i)}^2 \right\}$$
for some absolute constant $c > 0$.

\end{thm}

This is of course a generalization of \ref{4b}. In particular, the best upper bound on $\norm{v_k}^2$ one can hope to get by directly comparing it to $\norm{u_k}^2$ is approximately $k^{\frac{4}{3}}$.

The smallest example for \ref{4b} we construct in this paper is 14-dimensional, and it is easy to adjust it to get a 12-dimensional lattice by replacing $\frac{e_1 + e_6 + \cdots + e_{14}}{3}$ with $\frac{e_1 + e_4 + e_5 + \cdots + e_{12}}{3}$. It would be interesting to know what is the smallest dimension a lattice can have such that $\norm{v_n} > \overline{\lambda_n}$.

In an appendix which is joint with Lior Hadassi, we construct a lattice with a rather striking property.

\begin{thm}\label{Proj}

There exists a $43$-dimensional lattice $L$, such that any basis $b_1, \dots, b_{43}$ of $L$ which contains the shortest vector in the lattice is not the shortest basis, that is
$$\max \left\{ \norm{b_1}, \dots, \norm{b_{43}} \right\} > \overline{\lambda_{43}}.$$

\end{thm}

This is another generalization of Theorem \ref{4b}, because for this lattice $L$ not only does the "greedy algorithm" not produce the shortest basis, but any algorithm which at the first step takes the shortest vector in the lattice, does not produce the shortest basis. For example, in this lattice we have that
$$\max \left\{ \norm{u_1}, \dots, \norm{u_{43}} \right\} > \overline{\lambda_{43}}.$$

The lattice is constructed explicitly using a projective plane over a finite field, but the verification that it satisfies the above property was done with Sage.

\section{Proof of Theorem \ref{Bound}}

For the rest of this section, let $w_1, \dots, w_n \in L$ be a sequence of $n$ linearly independent vectors in $L$ such that $\norm{w_i} = \lambda_i$, and let $W_i = \mathrm{span}_{\mathbb{Z}} \left( w_1, \dots, w_i \right)$. For $W$ a subset of $\mathbb{R}^n$, let $P_W : \mathbb{R}^n \to \mathbb{R}^n$ denote the projection on $\mathrm{span}_{\mathbb{R}} \left( W \right)$. We will only use $P_W$ for lattices and linear subspaces. Recall that $v_1, \dots, v_n$ is the Minkowski basis of $L$, and let $V_i = \mathrm{span}_{\mathbb{Z}} \left( v_1, \dots, v_i \right)$.

First, we will prove a lemma which does not seem to be written explicitly in the literature, but the general principle has been used by van der Waerden \cite{Van Der Waerden} and Lagarias, Lenstra and Schnorr \cite{Bases} among others. To shorten our notation, let us say that $L'$ is a \emph{primitive sublattice} of $L$ if there exist linearly independent vectors $v_1, \dots, v_k \in L \setminus L'$ such that $\mathrm{span}_{\mathbb{Z}} \left( v_1, \dots, v_k, L' \right) = L$. Equivalently, one can complete a basis of $L'$ to a basis of $L$.

\begin{lem}\label{preHop}

Let $L_1$ be an $n$-dimensional lattice, and $L_2 \subset L_1$ a primitive $k$-dimensional sublattice. Let $y_1, \dots, y_k \in L_2$ be a sequence of $k$ linearly independent vectors, and let $b_1, \dots, b_k$ be their Gram-Schmidt orthogonalization. Let $B_i = \mathrm{span}_{R} \left( b_1, \dots, b_{i - 1}, b_{i + 1}, \dots, b_k \right)$ for $i = 1, \dots, k$. Then there exists a vector $y \in L$ such that $\mathrm{span}_{\mathbb{Z}} \left( L_2, y \right)$ is a primitive sublattice of $L_2$ such that
$$\norm{y}^2 \leq \max \left\{ \lambda_{k + 1}^2, \frac{\sum_{i = 1}^{k} \norm{P_{B_i} \left( y_i \right)}^2+ \lambda_{k + 1}^2}{4} \right\}.$$

\end{lem}

It is easier to understand the lemma in the case where $b_i = e_i$, or equivalently $y_i \in \mathrm{span}_{\mathbb{Z}} \left( e_1, \dots, e_i \right)$ which we can assume is the case after applying an orthogonal transformation. Then we can write for all $1 \leq i \leq k$,
$$y_i = \sum_{j = 1}^{i} y_{i, j} e_j$$
and our claim is that there exists a vector $y \in L$ such that $\mathrm{span}_{\mathbb{Z}} \left( L_2, y \right)$ is a primitive sublattice of $L_2$ of squared length at most
$$\max \left\{ \lambda_{k + 1}^2, \frac{\sum_{i = 1}^{k} y_{i, i}^2 + \lambda_{k + 1}^2}{4} \right\}.$$
We will prove this lemma, and then show some applications of it. Take a vector $y_0 \in L_1 \setminus L_2$ of length at most $\lambda_{k + 1}$. If $\mathrm{span}_{\mathbb{Z}} \left( L_2, y_0 \right)$ is a primitive sublattice of $L_1$, then $y = y_0$ and we are done. Otherwise, $P_{L_{2}^{\perp}} \left( y_0 \right)$ is not a primitive vector in $P_{L_{2}^{\perp}} \left( L_1 \right)$. Therefore, the shortest projection on $Y_k$ is of length at most $\frac{\norm{P_{L_{2}^{\perp}} \left( y_0 \right)}}{2}$.

Take a vector $y$ with a shortest nonzero projection on $L_{2}^{\perp}$. By subtracting an appropriate multiple of $y_k$ we can assume that the $e_k$ coordinate of $y$ is at most $\frac{\abs{y_{k, k}}}{2}$ in absolute value. Then, by subtracting the appropriate multiple of $y_{k - 1}$ we can assume that the $e_{k - 1}$ of $y$ is at most $\frac{\abs{y_{k - 1, k - 1}}}{2}$ and so on. Therefore, writing
$$y = \sum_{i = 1}^{k} \left[ y \right]_i e_i + w$$
where $w = P_{L_{2}^{\perp}} \left( y \right)$, we have
$$\norm{y}^2 = \sum_{i = 1}^{k} \left[ y \right]_{i}^2 + \norm{w}^2 \leq \sum_{i = 1}^{k} \frac{y_{i, i}^2}{4} + \frac{\norm{P_{L_{2}^{\perp}} \left( y_0 \right)}^2}{4}  \leq \frac{\sum_{i = 1}^{k} y_{i, i}^2 + \lambda_{k + 1}^2}{4}$$
which completes the proof of the lemma. \qed

Van der Waerden essentially used this lemma with $y_i = v_i$ for $i = 1, 2, \dots, k$ and bounded inductively $v_{i, i}^2 \leq \norm{v_i}^2 \leq \Delta_i \lambda_{i}^2$ to get
$$\norm{v_{k + 1}}^2 \leq \max \left\{ 1, \frac{\sum_{i = 1}^{k} \Delta_i + 1}{4} \right\} \lambda_{k + 1}^2.$$
Now it is a simply a matter of checking that $\Delta_k = \max \left\{ 1, \left( \frac{5}{4} \right)^{k - 4} \right\}$ satisfies the recurrence relation
$$\Delta_{k + 1} = \max \left\{ 1, \frac{\sum_{i = 1}^{k} \Delta_i + 1}{4} \right\}.$$
As we can see from this bound, any improvement on $\Delta_j$ for some $j$ gives us the kind of savings in Corollary \ref{Bound2}. Explicitly, suppose that for some $j$ we showed that
$$\Delta_j \leq \left( \frac{5}{4} \right)^{j - 4} - c$$
Then applying the above argument,
$$\Delta_{j + 1} \leq \left( \frac{5}{4} \right)^{j - 3} - \frac{c}{4}, \ \Delta_{j + 2} \leq \left( \frac{5}{4} \right)^{j - 2} - \frac{5 c}{16}, \dots$$
In general we get that
$$\Delta_{j + k} \leq \left( \frac{5}{4} \right)^{j + k - 4} - a_k$$
where $a_k$ satisfies the recurrence relation
$$a_0 = c, \ a_k = \frac{\sum_{i = 0}^{k - 1} a_i}{4}.$$
It is easy to check that for all $k > 0$ we have $a_k = c \cdot \frac{5^{k - 1}}{4^k}$ which means that we our bound is
$$\Delta_{j + k} \leq \left( 1 - c \cdot \frac{4^{j - 4}}{5^{j - 3}} \right) \cdot \left( \frac{5}{4} \right)^{j + k - 4}.$$
Substituting our improved bounds $\Delta_6 = \frac{6}{4}, \ \Delta_7 = \frac{7}{4}$ we find that
$$\Delta_k \leq \frac{608}{625} \cdot \left( \frac{5}{4} \right)^{k - 4}$$
as required. Therefore Theorem \ref{Bound} indeed implies Corollary \ref{Bound2}.

An important corollary of Lemma \ref{preHop} is the following:

\begin{cor}\label{Hop}

If $V_k = W_k$ then
$$\norm{v_{k + 1}}^2 \leq \frac{k + 1}{4} \lambda_{k + 1}^2.$$

\end{cor}

Note that this is exactly the conjectured bound. This corollary follows from Lemma \ref{preHop} by taking $y_i = w_i$ and using the inequality
$$\norm{\pi_{B_i \left( w_i \right)}} \leq \norm{w_i} = \lambda_i \leq \lambda_{k + 1}.$$
Now we will prove Theorem \ref{Bound}. By van der Waerden's bound, $\norm{v_i} \leq \lambda_i$ and therefore $v_i = w_i$ for $i l\leq 4$. We will first show the bound for $\norm{v_6}$. By Corollary \ref{Hop}, if $v_5 \in W_5$ then $V_5 = W_5$ and so $\norm{v_6} \leq \frac{6}{4} \lambda_{6}^2$. Therefore, we can assume that $v_5 \not\in W_5$. In particular, $v_1, v_2, v_3, v_4, w_5$ are not a primitive sequence of vectors, and therefore using the same argument as in Lemma \ref{preHop}, there is a vector $z_5 \in \mathrm{span}_{R} \left( w_1, \dots, w_5 \right)$ such that
$$\norm{P_{V_{4}^{\perp}} \left( z_5 \right)} \leq \frac{\norm{P_{V_{4}^{\perp}} \left( w_5 \right)}}{2} \leq \frac{\lambda_5}{2}, \ \norm{z_5}^2 \leq \frac{5}{4} \lambda_{5}^2$$
and $\mathrm{span}_{\mathbb{Z}} \left( V_4, z_5 \right)$ is a primitive sublattice of $L$. By our assumption, $z_5 \neq v_5$ and in fact $z_5 \not\in V_5$. If $\mathrm{span}_{\mathbb{Z}} \left( V_5, z_5 \right)$ is a primitive sublattice, then $\norm{v_6}^2 \leq \norm{z_5}^2 \leq \frac{5}{4} \lambda_{5}^2 < \frac{6}{4} \lambda_{6}^2$ as required. Otherwise, there exists a vector $z_6 \in \mathrm{span}_{\mathbb{Z}} \left( V_5, z_5 \right)$ such that $\mathrm{span}_{\mathbb{Z}} \left( V_5, z_6 \right)$ is a primitive sublattice of $L$ and
$$\norm{P_{V_{5}^{\perp}} \left( z_6 \right)} \leq \frac{\norm{P_{V_{5}^{\perp}} \left( z_5 \right)}}{2} \leq \frac{\norm{P_{V_{4}^{\perp}} \left( z_5 \right)}}{2} \leq \frac{\lambda_5}{4}.$$
Applying the same argument as in Lemma \ref{Hop} we can take $z_6$ to be of squared length at most
$$\norm{P_{V_{5}^{\perp}} \left( z_6 \right)}^2 + \frac{\sum_{i = 1}^{5} \norm{v_i}^2}{4} \leq \frac{11}{8} \lambda_{5}^2 < \frac{6}{4} \lambda_{6}^2.$$
as required.

Next, we will show the bound for $\norm{v_7}$. As before, if $V_6 = W_6$ then by Corollary \ref{Hop} we are done. Therefore we can assume that $v_6 \not\in W_6$. Now we divide into cases, depending on whether $V_5 = W_5$ or not.

Assume first that $V_5 = W_5$. As $v_6 \not\in W_6$, in particular $W_6$ is not a primitive sublattice, and therefore there exists a vector $z_6 \in \mathrm{span}_{\mathbb{Z}} \left( w_1, \dots, w_6 \right)$ such that
$$\norm{P_{V_{5}^{\perp}} \left( z_6 \right)} \leq \frac{\lambda_6}{2}, \ \norm{z_6}^2 \leq \frac{6}{4} \lambda_{6}^2$$
and $\mathrm{span}_{\mathbb{Z}} \left( V_5, z_6 \right)$ is a primitive sublattice. If $\mathrm{span}_{\mathbb{Z}} \left( V_6, z_6 \right)$ is a primitive sublattice, then $\norm{v_7}^2 \leq \norm{z_6}^2 \leq \frac{6}{4} \lambda_{6}^2 < \frac{7}{4} \lambda_{7}^2$ as required. Otherwise, as before there exists a vector $z_7$ such that
$$\norm{P_{V_{6}^{\perp}} \left( z_7 \right)} \leq \frac{\norm{P_{V_{5}^{\perp}} \left( z_6 \right) }}{2} \leq \frac{\lambda_6}{4}$$
such that $\mathrm{span}_{\mathbb{Z}} \left( V_6, z_7 \right)$ is a primitive sublattice of $L$. Using the same argument in Lemma \ref{Hop} where we take $y_i = w_i$ for $1 \leq i \leq 5$ and $y_6 = v_6$ we can assume that
\begin{align*}
    \norm{z_7}^2 \leq \frac{\sum_{i = 1}^{6} \norm{y_i}^2}{4} + \norm{P_{V_{6}^{\perp}} \left( z_7 \right)}^2 \leq \frac{4 + 4 + 4 + 4 + 4 + 6 + 1}{16} \lambda_{6}^2 = \frac{27}{16} \lambda_{6}^2\\
< \frac{7}{4} \lambda_{7}^2.
\end{align*}
Therefore, $\norm{v_7}^2 \leq \norm{z_7}^2 < \frac{7}{4} \lambda_{7}^2$

The case where $V_5 \neq W_5$ is similar, we get a vector $z_7$ such that
$$\norm{P_{V_{6}^{\perp}} \left( z_7 \right)} \leq \frac{\lambda_7}{4}$$
and $\mathrm{span}_{\mathbb{Z}} \left( V_6, z_7 \right)$ is a primitive sublattice of $L$. As above, we can assume that
$$\norm{v_7}^2 \leq \norm{z_7}^2 \leq \frac{4 + 4 + 4 + 4 + 5 + 6 + 1}{16} \lambda_{6}^2 = \frac{7}{4} \lambda_{6}^2 \leq \frac{7}{4} \lambda_{7}^2$$
which concludes the proof.

All that is left is to check when equality holds in the lemma. Whenever we have used Lemma \ref{preHop} we have always bounded the projection of $y_i$ by $\norm{y_i}$, and equality can only hold when $y_1, \dots, y_k$ are orthogonal. Furthermore, in order for the \ref{preHop} to be tight there needs to be a vector whose $e_i$ coordinate is exactly $\frac{y_{i, i}}{2}$. Finally, we have used $\lambda_i \leq \lambda_6, \lambda_7$ for the bounds on $v_6, v_7$ respectively.

It is easy to check that the equality case for our bound for $v_5, v_6$ can only occur when $V_5, V_6$ are similar to $D_{5}^{*}, D_{6}^{*}$ respectively. Looking at our bound for $v_7$, in the case where $V_6 = W_6$ it is clear that equality can only hold when $V_7$ is similar to $D_{7}^{*}$. The only other case that might give equaltiy is when $V_5 \neq W_5$ and $V_6 \neq W_6$, but there we have used both of the inequalities
$$\norm{v_5}^2 \leq \frac{5}{4} \lambda_{5}^2, \ \norm{v_6}^2 \leq \frac{6}{4} \lambda_{6}^2$$
and as we have just seen, equality cannot occur in these simultaneously. \qed

To get a significantly better bound for $\norm{v_k}^2$ for general $k$ a more delicate argument seems to be required. For example, in the case where $V_i \neq W_i$ for all $4 < i < k$, the improvement our argument gives as opposed to van der Waerden's bound is at best $\frac{\lambda_{k}^2}{4}$, which as we described above will not give a better bound than $\norm{v_k}^2 \leq c \cdot \left( \frac{5}{4} \right)^{k}$ for some constant $c$. Therefore, even
$$\norm{v_k}^2 = o \left( \left( \frac{5}{4} \right)^{k} \right)$$
would be a nontrivial improvement on our result.

\section{Proof of Theorems \ref{4b} and \ref{4b+}}

Let $p_1, p_2, \dots$ be the sequence of prime numbers in increasing order and define $a_{\ell} = p_{1}^2 + p_{2}^2 + \cdots + p_{\ell}^2 + 1$, and by convention $a_0 = 1$. Fix $k \geq 2$, and for each $1 \leq i \leq k$ denote
$$g_i = e_{a_{i - 1} + 1} + \cdots + e_{a_i}.$$
Define the lattices
$$L_k = \mathrm{span}_{\mathbb{Z}} \left( e_1, e_2, \dots, e_{a_k}, \frac{e_1 + g_1}{p_2}, \frac{e_1 + g_2}{p_2}, \dots, \frac{e_1 + g_k}{p_k} \right).$$
For example,
\begin{align*}
    L_2 = \mathrm{span}_{\mathbb{Z}} \left( e_1, \dots, e_{14}, \frac{e_1 + \cdots + e_5}{2}, \frac{e_1 + e_6 + e_7 + \cdots + e_{14}}{3} \right).
\end{align*}
As we will show in this section, these lattices will satisfies the conditions of Theorem \ref{4b+}. We will fix $k$, and compute the Minkowski and Korkin-Zolotarev bases for $L_k$. For convenience's sake, from now on let $L = L_k$ and $d = a_k = \dim (L)$.

Notice that for every $w \in L$ and every $1 \leq i \leq k$ there exist integers $0 \leq x_i (w) < p_i$ and a vector $z \in \mathbb{Z}^d$ such that
\begin{align*}
    w = z + \left( \sum_{i = 1}^{k} \frac{x_i (w)}{p_i} \right) e_1 + \sum_{i = 1}^{k} \frac{x_i (w)}{p_i} g_i
\end{align*}
because if we write
$$w = \sum_{i = 1}^{d} s_i(w) e_i + \sum_{i = 1}^{k} r_i(w) \frac{e_1 + g_i}{p_i}$$
then letting $x_i(w) = r_i(w) \ \bmod p_i$ works. Looking at the $p_i$-adic valuation of the $e_1$-coordinate, we see that if for some $i, \ x_i (w) \neq 0$ then the $e_1$ coordinate is not an integer. Therefore, if $w \in L$ is a non-integral vector (that is, $w \not\in \mathbb{Z}^n$) then the $e_1$ coordinate of $w$ is not an integer. This in a sense is the crux of the entire argument. Another fact that we see from this is that any non-integral vector $w \in L$ is of length strictly greater than 1, because if $x_i (w) \neq 0$ then writing $w = \sum_{j = 1}^{d} [w]_j e_j$ we have
\begin{align*}
    \norm{w}^2 = \sum_{j = 1}^{d} [w]_{j}^2 \geq [w]_{1}^2 + \sum_{j = a_{i - 1} + 1}^{a_i} [w]_{j}^2.
\end{align*}
For all $a_{i - 1} + 1 \leq j \leq a_i, \ [w]_{j}$ is a non-zero fraction with denominator $p_i$  and so $[w]_{j}^2 \geq \frac{1}{p_{i}^2}$. Similarly, $[w]_1$ is a non-zero fraction with denominator $p_1 p_2\cdots p_k$ and so $[w]_{1}^2 \geq \frac{1}{\prod_{i = 1}^{k} p_{i}^2}$. Combining this we get
\begin{align*}
    \norm{w}^2 \geq \frac{p_{i}^2}{p_{i}^2} + \frac{1}{\prod_{i = 1}^{k} p_{i}^2} > 1.
\end{align*}
Now we will compute the Minkowski basis. First we will prove that $e_2, \dots e_d$ is a primitive set of $d - 1$ vectors. Take a vector $w \in L \ \cap \ \mathrm{span}_{\mathbb{R}} \left( e_2, \dots, e_d \right)$. The $e_1$ coordinate of $w$ is 0 and in particular an integer, and therefore all coordinates of $w$ are integers which means that $w \in \mathrm{span}_{\mathbb{Z}} \left( e_2, \dots, e_d \right)$ as required. Furthermore as we have just shown these vectors are of the shortest length possible, and so the first $d - 1$ vectors in the Minkowski basis for $L$ will be
\begin{align*}
    v_i = e_{i + 1}
\end{align*}
for $1 \leq i < d$. Now let us see what the last vector $v_d$ can be. If for some $i$ we have $x_i \left( v_d \right) = 0$ then all the coordinates of every vector in $\mathrm{span}_{\mathbb{Z}} \left( v_1, \dots, v_d \right)$ will not have $p_i$ in the denominator, and so the $v_i$ cannot span $L$. Therefore, $x_i \left( v_d \right) \neq 0$ for all $1 \leq i \leq k$ and so writing $v_d = \sum_{i = 1}^{d} \left[ v_d \right]_{i} e_i$ we have
\begin{align*}
    \norm{v_d}^2 = \sum_{i = 1}^{d} \left[ v_d \right]_{i}^2 > \sum_{i = 1}^{k} \left( \sum_{j = a_{i - 1} + 1}^{a_i} \left[ v_d \right]_{j}^2 \right).
\end{align*}
Since for every $a_{i - 1} < j \leq a_i$ we have $\abs{\left[ v_{a_k} \right]_j} \geq \frac{1}{p_i}$,
\begin{align*}
    \norm{v_d}^2 > \sum_{i = 1}^{k} 1 = k.
\end{align*}
There is a basis of $L$ with shorter vectors: take for each $1 \leq i \leq k$ the vector
\begin{align*}
    \frac{e_1 + g_i}{p_i}
\end{align*}
and for all $1 \leq j \leq d$ such that $j \neq a_i$ for $0 \leq i < k$ the vector $e_j$. These are $d$ vectors in $L$, and all that we need to show in order to prove that this is a basis is that $e_1, e_{a_1}, e_{a_2}, \dots, e_{a_{k - 1}}$ are spanned by these vectors. This is true because
$$e_1 = p_k \frac{e_1 + e_{a_{k - 1} + 1} + \cdots + e_{a_k}}{p_k} - e_{a_{k - 1} + 1} - e_{a_{k - 1} + 2} - \cdots - e_{a_k}$$
and
\begin{align*}
    e_{a_j} = p_j \frac{e_1 + e_{a_{j - 1} + 1} + \cdots + e_{a_j}}{p_j} - e_{a_{j - 1} + 1} - \cdots - e_{a_j - 1} - e_1.
\end{align*}
Each $e_j$ is of length 1, and
\begin{align*}
    \norm{\frac{e_1 + e_{a_{i - 1} + 1} + \cdots + e_{a_i}}{p_i}}^2 = 1 + \frac{1}{p_{i}^2} \leq \frac{5}{4}
\end{align*}
which means that $\norm{v_{a_k}}^2 > \frac{4 k}{5} \overline{\lambda_{a_k}}^2$. In fact as we will see now, this is the Korkin-Zolotarev basis.

Now we will find the Korkin-Zolotarev basis. For convenience, denote $U_m = \mathrm{span}_{\mathbb{R}} \left( u_1, \dots, u_m \right)$. Like in the Minkowski basis, we have $u_1 = e_2$ and $u_2 = e_3$. However, $u_3 = \frac{e_1 + \cdots + e_5}{2}$. The orthogonal projections of $e_4, e_5$ on $U_{3}^{\perp}$ are $\frac{- e_1 + 2 e_4 - e_5}{3}, \frac{- e_1 - e_4 + 2 e_5}{3}$ respectively and therefore $u_4 = e_4$ and $u_5 = e_5$. I claim that the basis always looks like this, that is

\begin{lem}\label{basis}

For all $1 < j \leq k$ and for all $1 \leq \ell \leq p_{j}^2, \ \ell \neq 2$ we have
\begin{align*}
    u_{a_{j - 1} + \ell} = e_{a_{j - 1} + \ell}
\end{align*}
and
\begin{align*}
    u_{a_{j - 1} + 2} = \frac{e_1 + g_j}{p_j}.
\end{align*}

\end{lem}

From this lemma it immediately follows that
$$\max \left\{ \norm{u_1}^2, \dots, \norm{u_d}^2 \right\} = \frac{5}{4}.$$

We will prove the lemma by induction on $j$. Suppose we know the claim for all $j' < j$. As before, it is easy to see that $u_{a_{j - 1} + 1} = e_{a_{j - 1} + 1}$ and $u_{a_{j - 1} + 2} = \frac{e_1 + g_j}{p_j}$. Now let us prove by induction on $\ell$ that $u_{a_{j - 1} + \ell} = e_{a_{j - 1} + \ell}$. Assume that we have proved it for some $\ell$, and let us try to prove it for $\ell + 1$. We want to show that the shortest projection on $U_{a_{j - 1} + \ell}^{\perp}$ comes from $e_{a_{j - 1} + \ell + 1}$. Clearly the shortest projection on $U_{a_{j - 1} + \ell}^{\perp}$ will not come from a vector with a non-zero $e_q$ coordinate for $q > a_j$, and so it is sufficient to find the shortest vector that is spanned by the projections of $e_{a_{j - 1} + r}$. Notice that for $\ell < r \leq p_{j}^2$ the orthogonal projection of $e_{a_{j - 1} + r}$ on $U_{a_{j - 1} + \ell}^{\perp}$ is
\begin{align*}
    e_{a_{j - 1} + r} - \frac{e_{a_{j - 1} + \ell + 1} + \cdots + e_{a_j}}{p_{j}^2 - \ell}.
\end{align*}
Letting $m = p_{j}^2 - \ell$, renaming and re-scaling the vectors by $m$ what we want to show is that the shortest vector in the lattice
\begin{align*}
    \mathrm{span}_{\mathbb{Z}} \left( m e_1 - \left( e_1 + \cdots + e_m \right), m e_2 - \left( e_1 + \cdots + e_m \right), \dots, m e_m - \left( e_1 + \cdots + e_m \right) \right)
\end{align*}
is $m e_1 - \left( e_1 + \cdots + e_m \right)$. Every vector in this lattice is of the form
\begin{align*}
    w = \sum_{i = 1}^{m} a_i \left( m e_i - \left( e_1 + \cdots + e_m \right) \right) = \sum_{i = 1}^{m} \left( m a_i - \left( a_1 + \cdots + a_m \right) \right) e_i
\end{align*}
and its length squared is
\begin{align*}
    \norm{w}^2 = \sum_{i = 1}^{m} \left( m a_i - \left( a_1 + \cdots + a_m \right) \right)^2 = m^2 \left( \sum_{i = 1}^{m} a_{i}^2 \right) - m  \left( \sum_{i = 1}^{m} a_i \right)^2.
\end{align*}
We want to look for the $m$-tuples $(a_1, \dots, a_m) \in \mathbb{Z}^m$ which minimize this quantity. Notice that if there exist $i, j$ such that $a_i - a_j \geq 2$ then replacing $a_i$ with $a_i - 1$ and $a_j$ with $a_j + 1$ makes the length of $w$ smaller. Furthermore,
\begin{align*}
    \sum_{i = 1}^{m} \left( m e_i - \left( e_1 + \cdots + e_m \right) \right) = 0
\end{align*}
and therefore we can shift all the $a_i$ by a constant to assume that for all $i, \ a_i = 0, 1$. Suppose that there were $r$ 1's and $m - r$ 0's for some $0 < r < m$. Then,
\begin{align*}
    \norm{w}^2 = r m^2 - m r^2 = m r (m - r)
\end{align*}
which is clearly minimized by $r = 1, m - 1$. But if $r = m - 1$ then by shifting all the $a_i$ by 1 again, we get (up to sign) a vector of the form
\begin{align*}
    m e_q - \left( e_1 + \cdots + e_m \right)
\end{align*}
exactly as required, and so we have completed the proof of \ref{basis} \qed

As we have mentioned above, this shows that
$$\max \left\{ \norm{u_1}^2, \dots, \norm{u_d}^2 \right\} = \frac{5}{4}$$
and so
\begin{align*}
    \norm{v_{a_k}}^2 > \frac{4 k}{5} \norm{u_{a_k}}^2.
\end{align*}
All that is left is to understand the size of $a_k$ in terms of $k$. By the prime number theorem (though for our purposes a much weaker estimate is sufficient) $p_r \leq C r \log r$ for some constant $C > 0$ and so
\begin{align*}
    a_k = 1 + \sum_{r = 1}^{k} p_{r}^2 \leq \sum_{r = 1}^{k} p_{k}^2 \leq C^2 k^3 \log^2 k.
\end{align*}
Therefore for some constant $c > 0$
\begin{align*}
    k \geq \frac{c \cdot a_{k}^{\frac{1}{3}}}{\log^{\frac{2}{3}} a_k}
\end{align*}
and so
\begin{align*}
    \norm{v_{a_k}}^2 \geq \frac{c \cdot a_{k}^{\frac{1}{3}}}{\log^{\frac{2}{3}} a_k} \norm{u_{a_k}}^2
\end{align*}
as required.

\section{Appendix (joint with Lior Hadassi): Proof of Theorem \ref{Proj}}

First we will try to understand more about what the desired condition means. Let $w_1, \dots, w_n$ be one of the shortest bases of $L$, and let $w_0$ be the shortest vector in $L$. Then there exist some integers $a_1, \dots, a_n$ such that
$$w_0 = \sum_{i = 1}^{n} a_i w_i.$$
If $w_0$ does not participate in a shortest basis, then in particular it means that for all $i$ the set of vectors we get when we replace $w_i$ with $w_0$, that is
$$w_1, \dots, w_{i - 1}, w_0, w_{i + 1}, \dots, w_n$$
is not a basis of the lattice. As this set spans all the $w_j$ except for $w_i$, then the fact that it is not a basis is equivalent to saying that $w_i$ is not spanned by this set. Now, if $a_i = 0$ then this set of vectors is linearly dependent and so is not a basis. Otherwise,
$$w_i = \frac{w_0}{a_i} - \sum_{j \neq i}^{n} \frac{a_j}{a_i} w_i.$$
The right hand side is in $\mathrm{span}_{\mathbb{Z}} \left( w_1, \dots, w_{i - 1}, w_0, w_{i + 1}, \dots, w_n \right)$ if and only if $a_i = \pm 1$.

Suppose that $w_1, \dots, w_n$ was the shortest sequence of $n$ linearly independent vectors in $L$ such that no $w_i$ is a multiple of $w_0$, in the sense that it minimizes
$$\max \left\{ \norm{w_1}, \norm{w_2}, \dots, \norm{w_n} \right\}.$$
Then, the only other possible shortest bases of $L$ would be to replace one of the $w_i$ with $w_0$. Therefore in this case, if for all $i, \ a_i \neq \pm 1$ then $w_1, \dots, w_n$ is the unique shortest basis and in particular $w_0$ does not participate in any shortest basis.

Suppose we had some $n - 1$ dimensional lattice $L' \subset \mathbb{R}^{n - 1}$ that had exactly $2 n$ shortest vectors $\pm w_1, \dots, \pm w_{n}$, and that
$$L' = \mathrm{span}_{\mathbb{Z}} \left( w_1, \dots, w_n \right).$$
As $w_1, \dots, w_n$ are $n$ vectors in an $n - 1$ dimensional lattice, there exist coprime integers $a_1, \dots, a_n$ such that
$$\sum_{i = 1}^{n} a_i w_i = 0.$$
Take $\varepsilon_1, \dots, \varepsilon_n$ to be small non-zero irrational real numbers which are linearly independent over $\mathbb{Q}$. Our lattice will be
$$L = \mathrm{span}_{\mathbb{Z}} \left( w_1 + \varepsilon_1 e_n, w_2 + \varepsilon_2 e_n, \dots, w_n + \varepsilon_n e_n \right)$$
that is, we give a small "height" to each of $w_1, \dots, w_n$. The shortest vector in $L$ is clearly
$$\sum_{i = 1}^{n} a_i \left( w_i + \varepsilon_i e_n \right) = \left( \sum_{i = 1}^{n} a_i \varepsilon_i \right) e_n$$
as the projection on $\mathbb{R}^{n - 1}$ of any other vector which is not a multiple of it is significantly longer. The same argument, that is looking at the lengths of projections of vectors in $L$, shows that
$$w_1 + \varepsilon_1 e_n, \ w_2 + \varepsilon_2 e_n, \dots, w_n + \varepsilon_n e_n$$
is the unique shortest sequence of $n$ linearly independent vectors in $L$ such that none of them are a multiple of the shortest vector in $L$. Therefore, if none of the $a_i$ are equal to $\pm 1, L$ satisfies the condition that we want.

We will now show two attempts at finding a lattice $L'$ which satisfies the above conditions. The first attempt does not work, but we have included it here in order to gain some intuition.

In both of these constructions we use a finite projective plane, and the idea behind this is that the highly symmetrical form of the projective plane should help prevent additional shorter vectors: if there was a short vector in our lattice which was not one of the $n + 1$ vectors we used to define it, by symmetry we would have many short vectors, and this would lead to a contradiction. "By chance" the first, simpler example doe not work, and so we must choose a slightly more complicated lattice which does work.

We will take $L' \subset \mathbb{Z}^{21}$ to be a 21-dimensional lattice. All of our vectors will be of the form $e_i + e_j + e_k$, and so we just need to specify the triplets $(i, j, k)$. Take three copies of $\mathbb{P}^{2} \left( \mathbb{F}_2 \right)$, the projective plane over the field of two elements, and label the points of the first plane $1, 2, \dots, 7$, the points of the second plane $8, \dots, 14$ and the points of the third plane $15, \dots, 21$. Each line in $\mathbb{P}^{2} \left( \mathbb{F}_2 \right)$ contains $3$ points, and by duality the number of points and lines is the same, that is $7$. So taking all the lines contained in one of the three copies gives us $21$ triplets. Finally, take the triplet $(1, 8, 15)$.

For convenience, let $L_{proj}$ denote the $\mathrm{span}_{\mathbb{Z}}$ in $\mathbb{R}^7$ of the $7$ vectors corresponding to the $7$ lines in $\mathbb{P}^{2} \left( \mathbb{F}_2 \right)$.

We have chosen $22$ vectors in $\mathbb{Z}^{21}$, and we let $L'$ be the $\mathrm{span}_{\mathbb{Z}}$ of these vectors. First, let us show that $L'$ is of full rank in $\mathbb{R}^{21}$. For this it is enough to show that the $L_{proj}$ is of full rank, or equivalently $\mathrm{covol} \left( L_{proj} \right) \neq 0$. We can do this by computing the Gram matrix: every line contains $3$ points, and every two distinct lines intersect in one point. Therefore the Gram matrix is a $7 \cross 7$ matrix with $3$'s on the diagonal and $1$'s off the diagonal. It is easy to see that the eigenvalues of this matrix are $2$ with multiplicity $6$ and $9$ with multiplicity $1$, and so the determinant of the Gram matrix is $2^6 \cdot 3^2$. As the determinant of the Gram matrix is the square of $\mathrm{covol} \left( L_{proj} \right)$, we have $\mathrm{covol} \left( L_{proj} \right) = 24$ and in particular it is nonzero.

We will state two useful facts about $L_{proj}$. The first is that
$$L_{proj} \subset \left\{ \left( x_1, x_2, \dots, x_7 \right) \in \mathbb{Z}^n \ : \ x_1 + x_2 + \cdots + x_7 = 0 \ (\bmod \ 3) \right\}.$$
This is clear, because each of the basis vectors of $L_{proj}$ are in the right hand side, and the condition is preserved under integral linear combinations of vectors.

The second fact is that for all $1 \leq i \neq j \leq 7$ we have $e_i - e_j \not \in L_{proj}$. Suppose to the contrary that some $e_{i_0} - e_{j_0} \in L_{proj}$. Then, as projective transformations on $\mathbb{P}^{2} \left( \mathbb{F}_2 \right)$ are transitive on pairs of points we would have $e_i - e_j \in L_{proj}$ for all $1 \leq i \neq j \leq 7$. This and the first fact together imply that
$$L_{proj} = \left\{ \left( x_1, x_2, \dots, x_7 \right) \in \mathbb{Z}^7 \ : \ x_1 + x_2 + \cdots + x_7 = 0 \ (\bmod \ 3) \right\}$$
but we know that their volumes are different: $\mathrm{covol} \left( L_{proj} \right) = 24$ whereas
$$\mathrm{covol} \left( \left\{ \left( x_1, x_2, \dots, x_7 \right) \in \mathbb{Z}^7 \ : \ x_1 + x_2 + \cdots + x_7 = 0 \ (\bmod \ 3) \right\} \right) = 3$$
and this is a contradiction.

Now we will show that the $22$ vectors we have chosen are indeed the $22$ shortest vectors in $L'$. Take a vector $w \in L'$ of length at most $\sqrt{3}$. As $w \in \mathbb{Z}^n$, all the nonzero coordinates of $w$ are equal to $\pm 1$, and there are at most $3$ of those. Furthermore, the sum of the coordinates of $w$ is divisible by $3$, as it is a linear combination of vectors whose sum of coordinates is divisible by $3$, and therefore either
$$w = e_i - e_j$$
or up to sign
$$w = e_i + e_j + e_k.$$
Write $w = w_1 + w_2 + w_3$, where
$$w_i \in \mathrm{span}_{\mathbb{Z}} \left( e_{7 i - 6}, \dots, e_{7 i} \right)$$
that is, we split $w$ into its components which lie in each of the 3 copies of the projective plane. Notice that the sum of the coordinates of all the $w_i$ are equal, as $w \equiv c \left( e_1 + e_{8} + e_{15} \right) \ \left( \bmod \ 3 \right)$ for some constant $c \in \mathbb{Z}$. If $w = e_i - e_j$ this implies that $i, j$ are in the same projective. Writing $w = m \left( e_1 + e_8 + e_{15} \right) + w'$ for some $m \in \mathbb{Z}$, where $w'$ is a linear combination of vectors corresponding to lines in one of the projective planes, then we see that
$$m e_1 \in L_{proj}$$
and
$$m e_1 + e_i - e_j \in L_{proj}$$
which means that
$$e_i - e_j \in L_{proj}$$
which is a contradiction. If $w = e_i + e_j + e_k$, then either all $i, j, k$ come from the same projective plane or each of them come from a different one. If all $i, j, k$ are from the same projective plane, then letting $l$ be the third point on the line through $i, j$ we have $e_i + e_j + e_l \in L_{proj}$ which means that $e_k - e_l \in L'$, which we have already shown is a contradiction. The final case is taken care of in the same way.
\\
Let our $22$ vectors be $v_1, \dots, v_{22}$. If we take the coprime integers $a_1, \dots, a_{22}$ such that
$$\sum_{i = 1}^{22} a_i v_i = 0$$
then all we need in order for $L'$ to be an example is that none of the $a_i$ are equal to $\pm 1$. Unfortunately, this is not the case. More specifically, it turns out that the coefficient of $e_1 + e_8 + e_{15}$ is $-6$, the coefficient of any vector which corresponds to a line containing one of the points $1, 8, 15$ is $2$, and the coefficient of the rest of the vectors is $- 1$.

Now we come to our real example. $L' \subset \mathbb{Z}^{42}$ is a $42$-dimensional lattice, which we define as the span of $43$ vectors. Each of this vectors is of the form $e_{i_1} + e_{i_2} + e_{i_3} + e_{i_4} + e_{i_5}$, and so we just need to specify the quintuplets $(i_1, i_2, i_3, i_4, i_5)$. Similarly to our discussion above, take two copies of $\mathbb{P}^{2} \left( \mathbb{F}_4 \right)$, the projective plane over the field of $4$ elements. $42$ of our quintuplets will consist of the $42$ lines in these two copies of $\mathbb{P}^{2} \left( \mathbb{F}_4 \right)$, and the last quintuplet consists of $3$ points in the first copy of the projective plane which are not colinear, and $2$ points in the second copy of the projective plane.

In principle it is possible to check by hand that this lattice works in the same way as above, but the number of cases becomes very large. This was verified using Sage, and the code is written below. The algorithm is as follows:

We define our 43 vectors $b_1, \dots, b_{43}$. Then, we compute the linear dependence between these vectors and check that none of the coefficients are $\pm 1$. After that, we need to check that $b_1, \dots, b_{43}$ are the only shortest vectors in the lattice. The lattice they span is a subset of
$$\left\{ (x_1, x_2, x_3, x_4, x_5) \in \mathbb{Z}^5 : x_1 + x_2 + x_3 + x_4 + x_5 \equiv 0 ( \bmod 5) \right\}$$
and therefore any shorter vector must be one of $e_i - e_j, \ e_i + e_j - e_k - e_m, \ e_i + e_j + e_k + e_m + e_n$ for some indices $i, j, k, m, n$. The code runs over all such vectors, and for each vector it checks if it is an integral linear combination of $b_1, \dots, b_{43}$. To do this, given a vector $v$ we first express it in a unique way as a rational linear combination of $b_2, \dots, b_{43}$. This can be done efficiently by a one time computation of the inverse matrix of $b_2, \dots, b_{43}$, which we call $M^{- 1}$, and then the vector $v \cdot M^{- 1}$ is the desired linear combination. Now, the linear dependence is of the form
$$\sum_{i = 1}^{43} a_i b_i = 0$$
or equivalently
$$b_1 + \sum_{i = 2}^{43} \frac{a_i}{a_1} b_i = 0$$
so if $v$ is an integral linear combination of $b_1, \dots, b_{43}$, this combination must be $v \cdot M^{- 1}$ plus $k$ times the relation $b_1 + \sum_{i = 2}^{43} \frac{a_i}{a_1}$ for some $0 \leq k < a_1$, and so we check these $a_1$ combinations and verify that all of them contain a non-integral coordinate.
\\
The code is far from being optimized. For example, just replacing $b_2, \dots, b_{43}$ with $b_1, b_3, \dots, b_{43}$ would shorten the runtime by a factor of approximately $\frac{3}{2}$, as $a_1 = 3$ and $a_2 = 2$. Furthermore, we have not at all utilized the high symmetry of the projective plane. As a ballpark estimate, all triplets of noncolinear points are projectively equivalent to each other, which should reduce the number of cases by a factor of at least $\binom{21}{3} = 1330$, and combining the computer calculations with arguments as we did above could give an even better improvement. In our case this did not matter as the code took 13 minutes to run, but for further investigations optimizations like this may be worthwhile.

\begin{verbatim}
#For vector v and list w of length N, checks if some v + j*w is integral,
#for maxrange > j >= 0
def integerrel(v, w, maxrange, N):
    wvector = vector(w)
    for j in range(maxrange):
        if (v + j*wvector in ZZ^N):
            return True
    return False

#Checks if the lattice spanned by the N+1 vectors in the list l satisfies our conditions   
def check(l, N):
    #Initializing Q^N and its standard basis, converting l to a list of vectors
    basis = []
    for i in range(N):
        vec = vector(ZZ, N)
        vec[i] = 1
        basis.append(vec)

    veclist = []
    for j in l:
        vec = vector(ZZ, N)
        for i in j:
            vec += basis[i]
        veclist.append(vec)

    V = QQ ^ N
    
    #Checking if the linear relation of our vectors contains plus or minus 1
    relations = V.linear_dependence(veclist)
    if len(relations) != 1:  # Our 43 vectors should span a 42-dimensional space,
            # so there should be one linear dependence up to a constant
        return "Lattice is not of full rank"
    relation = [i for i in relations[0]]
    for j in relation:
        if j == 1 or j == -1:
            return "Linear combination contains plus or minus one"
    print(relation)  # Just in case
    
    #Taking the first coefficient in our relation and dividing by it
    #Make sure that it is not 0! In our case all coefficients are nonzero
    first_coef = relation[0]
    del relation[0]
    for j in range(len(relation)):
        relation[j] = relation[j] / first_coef
    
    #Taking the inverse matrix of the N x N matrix which corresponds to N of our vectors    
    Mat = matrix(QQ, veclist[1:])
    Inv = Mat.solve_right(matrix.identity(N))
    
    #Multiplying a matrix and a vector is costly. To compute (e_i + e_j - e_k - e_m)*Inv,
    #We just need to precompute e_i*Inv and add vectors, which takes less time
    invbasis = []
    for i in range(N):
        v = vector(ZZ, N)
        v[i] = 1
        invbasis.append(v * Inv)
    
    #Checking that there are no shorter vectors in the lattice
    for i in range(N):
        for j in range(i + 1, N):
            v = invbasis[i] - invbasis[j]
            if integerrel(v, relation, first_coef, N):
                return "Lattice contains a vector of the form e_i - e_j"

    for i in range(N):
        for j in range(i + 1, N):
            for k in range(j + 1, N):
                for m in range(k + 1, N):
                    v = invbasis[i] + invbasis[j] - invbasis[k] - invbasis[m]
                    if integerrel(v, relation, first_coef, N):
                        return "Lattice contains a vector", \
                                "of the form e_i + e_j - e_k - e_m"
                    v = invbasis[i] - invbasis[j] + invbasis[k] - invbasis[m]
                    if integerrel(v, relation, first_coef, N):
                        return "Lattice contains a vector", \
                               "of the form e_i - e_j + e_k - e_m"
                    v = invbasis[i] - invbasis[j] - invbasis[k] + invbasis[m]
                    if integerrel(v, relation, first_coef, N):
                        return "Lattice contains a vector", \
                               "of the form e_i - e_j - e_k + e_m"

    for i in range(N):
        for j in range(i + 1, N):
            for k in range(j + 1, N):
                for m in range(k + 1, N):
                    for n in range(m + 1, N):
                        if (i, j, k, m, n) not in l:
                            v = invbasis[i] + invbasis[j] + invbasis[k] + \
                                invbasis[m] + invbasis[n]
                            if integerrel(v, relation, first_coef, N):
                                return "Lattice contains a vector", \
                                       "of the form e_i + e_j + e_k + e_m + e_n"
    return "Yay!"    

def main():
    #Initialize points in the projective plane
    Points = list(ProjectiveSpace(2)/GF(4))
    Lines = []
    for i in range(len(Points)):
        #By duality, each line is the set of points orthogonal to a fixed point i
        #For each point i, tup is the dual line
        #tup2 is the dual line in the second copy of the projective plane
        tup = tuple(j for j in range(len(Points)) if (Points[j][0]*Points[i][0] + \ 
                            Points[j][1]*Points[i][1] + Points[j][2]*Points[i][2]) == 0)
        tup2 = tuple(j + 21 for j in tup)
        Lines.append(tup)
        Lines.append(tup2)
    Lines.append((0, 1, 4, 21, 22))
    N = 42
    return check(Lines, N)

if __name__ == "__main__":
    main()
\end{verbatim}
\textbf{Acknowledgements}: I thank Barak Weiss for his course on Geometry of Numbers and for his many comments on the numerous drafts of this paper, both of which are directly responsible for this paper being written. I also thank Achill Sch\"{u}rmann for some helpful comments and encouragement, and Eyal Litvin for suggesting to me the problem which became Theorem \ref{Proj}. I thank Jacques Martinet for some helpful comments. Finally, I thank many of my friends for helping me cope with various computer-related issues.

\end{document}